\documentclass[leqno]{article}

\begin{document}

\title{Some algebraic notions related to analysis on metric spaces}

\author{Stephen Semmes	\\
	Rice University}

\date{}

\maketitle

\begin{abstract}
There are versions of ``calculus'' in many settings, with various
mixtures of algebra and analysis.  In these informal notes we consider
a few examples that suggest a lot of interesting questions.
\end{abstract}

\tableofcontents

\section{Functions on spaces}
\label{functions, spaces}
\setcounter{equation}{0}

	Let $X$ be a topological space, and let $\mathcal{C}(X)$ be
the algebra of continuous real-valued functions on $X$.  For each $p
\in X$, the collection $\mathcal{I}_p(X)$ of $f \in \mathcal{C}(X)$ such
that $f(p) = 0$ is an ideal in $\mathcal{C}(X)$, since the product of
two functions $f_1, f_2 \in \mathcal{C}(X)$ vanishes at $p$ if $f_1(p)
= 0$ or $f_2(p) = 0$.

	If $X$ is the $n$-dimensional Euclidean space ${\bf R}^n$ and
$f$ is a continuously-differentiable function which vanishes at $p$,
then
\begin{equation}
\label{f(x) = O(|x - p|)}
	f(x) = O(|x - p|)
\end{equation}
for $x$ near $p$.  A similar statement holds for functions on smooth
manifolds.

	In general, a function $f \in \mathcal{I}_p(X)$ does not have
to vanish at $p$ at any particular rate.  One might specify a rate at
which a function vanishes with a condition like (\ref{f(x) = O(|x -
p|)}).  On a metric space, one might work with functions in a
Lipschitz class, for instance, which implies such a condition.

	If $f$ is a polynomial on ${\bf R}^n$ such that $f(p) = 0$,
then there are polynomials $f_1, \ldots, f_n$ which satisfy
\begin{equation}
	f(x) = (x_1 - p_1) \, f_1(x) + \cdots + (x_n - p_n) \, f_n(x),
\end{equation}
where $x_1, \ldots, x_n$ are the standard coordinates of $x \in {\bf
R}^n$.  Of course, this implies (\ref{f(x) = O(|x - p|)}).

	This simple algebraic statement works just as well for
polynomials with coefficients in any field.  One can also consider
rational functions where the denominators are nonzero.

	If $f$ is a continuously-differentiable function on ${\bf R}^n$,
then
\begin{equation}
\label{f(x) = o(|x - p|)}
	f(x) = o(|x - p|)
\end{equation}
if and only if $f(p) = 0$ and $df_p = 0$, where $df_p$ denotes the
differential of $f$ at $p$.  If $f$ is continuously-differentiable of
order two, then it follows that
\begin{equation}
\label{f(x) = O(|x - p|^2)}
	f(x) = O(|x - p|^2).
\end{equation}

	For any continuously-differentiable function $f$ on ${\bf R}^n$,
\begin{equation}
	f(x) - f(p) - df_p(x - p) = o(|x - p|),
\end{equation}
and
\begin{equation}
	f(x) - f(p) - df_p(x - p) = O(|x - p|^2)
\end{equation}
when $f$ is continuously differentiable of order two.  Thus the
difference between vanishing to first order at $p$ and vanishing more
quickly is exactly described by the differential of $f$ at $p$.

	For a polynomial $f$ in $n$ variables, the idea that $f$
vanishes at $p$ to second order can be expressed algebraically by
\begin{equation}
	f(x) = \sum_{j, l = 1}^n (x_j - p_j) \, (x_l - p_l) \, f_{j, l}(x),
\end{equation}
where $f_{j, l}$ are also polynomials.

	On any metric space, it makes sense to talk about a function
vanishing to first order, more quickly than first order, or second
order at a point, as in (\ref{f(x) = O(|x - p|)}), (\ref{f(x) = o(|x -
p|)}), and (\ref{f(x) = O(|x - p|^2)}).  Depending on the
circumstances, there may be nice classes of regular functions with
some kind of derivatives, suitable versions of Taylor's theorem, and
so on.

	The real and complex numbers are equipped with metrics for
which algebraic and analytic notions of vanishing to some order are
compatible.  There is a similar metric on the $p$-adic numbers.

	There are translation-invariant metrics on ${\bf R}^n$
associated to nonisotropic dilations which lead to different
measurements of vanishing, smoothness, and degrees of polynomials.
This can be extended further to left-invariant metrics on nilpotent
Lie groups.  Just as a smooth manifold is approximately like a
Euclidean space locally, and thus has a version of calculus, there are
versions of calculus for sub-Riemannian spaces connected to those for
nilpotent Lie groups.

\section{Exponentiation}
\label{exponentiation}
\setcounter{equation}{0}

	The classical exponential function on the real line can be
defined by the power series
\begin{equation}
	\exp x = \sum_{n = 0}^\infty \frac{x^n}{n!}.
\end{equation}
It can also be characterized by the differential equation $f' = f$
with the condition $f(0) = 1$.

	This makes sense as a formal power series with coefficients in
any field of characteristic $0$ too.  The radius of convergence on the
$p$-adic numbers is positive and can be determined explicitly.

	On an open set in ${\bf R}^n$, or any smooth manifold, one can
consider a vector field $V$ with smooth coefficients.  One can view
this as a derivation on the algebra of smooth functions on the
underlying open set or manifold.

	Well-known results about systems of ordinary differential
equations imply the local existence and uniqueness of integral curves
for $V$.  Under suitable conditions, one has global solutions.

	One can also look at this in terms of the partial differential
equation
\begin{equation}
	\frac{\partial}{\partial t} f = V \, f,
\end{equation}
where $f(x, t)$ is a smooth function of $x$ in the open set or
manifold and $t \in {\bf R}$.

	If there are global solutions, as in the case of a compact
manifold without boundary, then one is basically exponentiating $V$ as
a derivation on the algebra of smooth functions on the underlying
space to get a one-parameter group of automorphisms of the algebra.
These automorphisms can be described by composing smooth functions on
the space with diffeomorphisms, where the diffeomorphisms correspond
to the flow determined by the vector field.

	A vector field on a sub-Riemannian space may be considered as
a vector field on an ordinary underlying manifold.  However, there are
additional conditions for compatibility with the sub-Riemannian
structure.

	Integrating vector fields on other metric spaces can be a
fascinating enterprise.

	For any linear transformation $L$ on a vector space $A$, one
can try to exponentiate $L$ or scalar multiples of $L$ to get
invertible linear transformations on $A$.  If $A$ is an algebra and
$L$ is a derivation, then the exponential ought to be an automorphism
of $A$.

	As in the case of the classical exponential function, one can
try to approach exponentiation of linear transformations in terms of
power series or differential equations.  If $A$ is a vector space over
a field of characteristic $0$, then
\begin{equation}
	\exp (t \, L) = \sum_{n = 0}^\infty \frac{t^n \, L^n}{n!}
\end{equation}
makes sense as a formal power series in $t$ with coefficients in
linear transformations on $A$.  If $A$ is a Banach space and $L$ is a
bounded linear transformation on $A$, then this series converges
absolutely for all $t$ in the Banach space of bounded linear
transformations on $A$.

	One might think of a vector field as a continuous linear
transformation on a topological vector space, or as an unbounded
linear transformation on a subspace of a Banach space.  In this case,
one normally approaches the exponential in terms of differential
equations.

\section{Buildings and symmetric spaces}
\label{buildings and symmetric spaces}
\setcounter{equation}{0}

	Classical symmetric spaces are homogeneous spaces obtained
from quotients of semi-simple Lie groups in a special way.  They are
quite interesting in particular for actions on them by discrete
subgroups of semisimple Lie groups.

	Buildings are generalizations of symmetric spaces with nice
properties for actions by broader classes of discrete groups.  Some of
these discrete groups are related to $p$-adic matrix groups in much the
same way as for real Lie groups in the classical case.

	In the classical situation, the boundary of a rank $1$
symmetric space of negative curvature is related to ordinary Euclidean
geometry and the geometry of sub-Riemannian spaces in a natural way.
In \cite{b-p-1}, it is shown that the boundaries of certain buildings
are fractal spaces which enjoy Poincar\'e inequalities similar to
those on Euclidean and sub-Riemannian spaces.  These fractal boundary
spaces have topological dimension $1$ and are far from being manifolds
even topologically.

	Of course, one can look at versions of calculus on discrete
structures too, on discrete groups in particular.

\section{Locally compact groups}
\label{locally compact groups}
\setcounter{equation}{0}

	Locally compact topological groups are very attractive as
general objects on which to consider basic notions of harmonic
analysis.  In particular, one always has Haar measure and
convolutions.  However, the solution of Hilbert's fifth problem shows
that such a group is a Lie group under some conditions related to
connectivity and finite-dimensionality.  Thus, there is so much
structure that something relatively close to the classical situation
may very well already be in the classical situation.

\section{Holomorphic functions}
\label{holomorphic functions}
\setcounter{equation}{0}

	A basic advantage of holomorphic functions as compared to
harmonic functions is that the product of two holomorphic functions is
holomorphic.  This leads to algebras of holomorphic functions, modules
over such algebras, and so on.  Local versions of this are closely
connected to formal power series.

	One can also look at holomorphic functions on suitable classes
of complex analytic metric spaces.  Let us suppose that there is some
local regularity for holomorphic functions.  It seems natural to ask
about analogues of power series expansions or related properties.
Some information may follow from appropriate elliptic analysis.

\section{Quasiconformal mappings}
\label{quasiconformal mappings}
\setcounter{equation}{0}

	In the theory of quasiconformal mappings in the plane, there
are remarkable existence results for mappings with prescribed
dilatation.  These results imply that perturbations of the classical
Cauchy--Riemann equations are in a sense equivalent to the original
ones.  Specifically, holomorphic functions with respect to perturbed
equations can be expressed as compositions of ordinary holomorphic
functions with quasiconformal mappings under suitable conditions.
This probably does not work very well in general for complex analytic
metric spaces, even when the complex dimension is equal to $1$.

\section{Uniform algebras}
\label{uniform algebras}
\setcounter{equation}{0}

	It seems natural to take a fresh look at deformation theory of
uniform algebras, as in \cite{ja1, ja2, ja3, ja4, ja5, ro1, ro2, ro3,
ro4}, in connection with complex analytic metric spaces.


\begin{thebibliography}{80}


\bibitem {a} L.~Ahlfors, {\it Lectures on Quasiconformal Mappings},
2nd edition, American Mathematical Society, 2006.

\bibitem {a-s} L.~Ambrosio and F.~Serra Cassano, editors, {\it
Lectures Notes on Analysis in Metric Spaces}, Scuola Normale
Superiore, Pisa, 2000.

\bibitem {a-t} L.~Ambrosio and P.~Tilli, {\it Topics on Analysis in
Metric Spaces}, Oxford University Press, 2004.

\bibitem {a-m} M.~Atiyah and I.~Macdonald, {\it Introduction to
Commutative Algebra}, Addison-Wesley, 1969.

\bibitem {a-c-g} P.~Auscher, T.~Coulhon, and A.~Grigoryan, editors,
{\it Heat Kernels and Analysis on Manifolds, Graphs, and Metric
Spaces}, American Mathematical Society, 2003.

\bibitem {b-g} R.~Beals and P.~Greiner, {\it Calculus on Heisenberg
Manifolds}, Princeton University Press, 1988.

\bibitem {b-r} A.~Bella\"{\i}che and J.-J.~Risler, editors, {\it
Sub-Riemannian Geometry}, Birkh\"auser, 1996.

\bibitem {bir-m} G.~Birkhoff and S.~Mac Lane, {\it A Survey of Modern
Algebra}, 4th edition, Macmillan, 1977.

\bibitem {boc-m} S.~Bochner and W.~Martin, {\it Several Complex
Variables}, Princeton University Press, 1948.

\bibitem {b} M.~Bourdon, {\it Immeubles hyperboliques, dimension
conforme et rigidit\'e de Mostow}, Geometric and Functional Analysis
{\bf 7} (1997), 245--268.

\bibitem {b-p-1} M.~Bourdon and H.~Pajot, {\it Poincar\'e inequalities
and quasiconformal mappings on the boundary of some hyperbolic
buildings}, Proceedings of the American Mathematical Society {\bf 127}
(1999), 2315--2324.

\bibitem {b-p-2} M.~Bourdon and H.~Pajot, {\it Cohomologie $l_p$ et
espaces de Besov}, Journal f\"ur die Reine und Angewandte Mathematik
{\bf 558} (2003), 85--108.

\bibitem {br} A.~Browder, {\it Introduction to Function Algebras},
Benjamin, 1969.

\bibitem {bn} K.~Brown, {\it Buildings}, Springer-Verlag, 1998.

\bibitem {cr} H.~Cartan, {\it Elementary Theory of Analytic Functions
of One or Several Variables}, Dover, 1995.

\bibitem {cs} J.~Cassels, {\it Local Fields}, Cambridge University
Press, 1986.

\bibitem {ch} J.~Cheeger, {\it Differentiability of Lipschitz
functions on metric measure spaces}, Geometry and Functional Analysis
{\bf 9} (1999), 428--517.

\bibitem {c-w-1} R.~Coifman and G.~Weiss, {\it Analyse Harmonique
Non-Commutative sur Certains Espaces Homog\`enes}, Lecture Notes in
Mathematics {\bf 242}, Springer-Verlag, 1971.

\bibitem {c-w-2} R.~Coifman and G.~Weiss, {\it Extensions of Hardy
spaces and their use in analysis}, Bulletin of the American
Mathematical Society {\bf 83} (1977), 569--645.

\bibitem {cns} A.~Connes, {\it Noncommutative Geometry}, Academic
Press, 1994.

\bibitem {e} D.~Eisenbud, {\it Commutative Algebra: With a View toward
Algebraic Geometry}, Springer-Verlag, 1995.

\bibitem {fal} K.~Falconer, {\it The Geometry of Fractal Sets},
Cambridge University Press, 1986.

\bibitem {fed} H.~Federer, {\it Geometric Measure Theory},
Springer-Verlag, 1969.

\bibitem {f-s} G.~Folland and E.~Stein, {\it Hardy Spaces on
Homogeneous Groups}, Princeton University Press, 1982.

\bibitem {ga1} T.~Gamelin, {\it Uniform Algebras}, Prentice-Hall,
1969.

\bibitem {ga2} T.~Gamelin, {\it Uniform Algebras and Jensen Measures},
Cambridge University Press, 1978.

\bibitem {g-j} L.~Gillman and M.~Jerison, {\it Rings of Continuous
Functions}, Springer-Verlag, 1976.

\bibitem {g} F.~Gouv\^ea, {\it $p$-Adic Numbers: An Introduction}, 2nd
edition, Springer-Verlag, 1997.

\bibitem {g-r} H.~Grauert and R.~Remmert, {\it Coherent Analytic
Sheaves}, Springer-Verlag, 1984.

\bibitem {gng} R.~Gunning, {\it Introduction to Holomorphic Functions
of Several Complex Variables}, volumes I, II, III, Wadsworth \& Brooks
/ Cole, 1990.

\bibitem {gng-r} R.~Gunning and H.~Rossi, {\it Analytic Functions of
Several Complex Variables}, Prentice-Hall, 1965.

\bibitem {h-k} P.~Haj{\l}asz and P.~Koskela, {\it Sobolev Met
Poincar\'e}, Memoirs of the American Mathematical Society {\bf 688},
2000.

\bibitem {h1} J.~Heinonen, {\it Lectures on Analysis on Metric
Spaces}, Springer-Verlag, 2001.

\bibitem {h2} J.~Heinonen, {\it Lectures on Lipschitz analysis},
Reports of the Department of Mathematics and Statistics {\bf 100},
University of Jyv\"askyl\"a, 2005.

\bibitem {h3} J.~Heinonen, {\it Nonsmooth calculus}, to appear,
Bulletin of the American Mathematical Society.

\bibitem {h} L.~H\"ormander, {\it An Introduction to Complex Analysis
in Several Variables}, 3rd edition, North-Holland, 1990.

\bibitem {ja1} K.~Jarosz, {\it The uniqueness of multiplication in
function algebras}, Proceedings of the American Mathematical Society
{\bf 89} (1983), 249--253.

\bibitem {ja2} K.~Jarosz, {\it Perturbations of uniform algebras},
Bulletin of the London Mathematical Society {\bf 15} (1983), 133--138.

\bibitem {ja3} K.~Jarosz, {\it Perturbations of uniform algebras II},
Journal of the London Mathematical Society (2) {\bf 31} (1985),
555-560.

\bibitem {ja4} K.~Jarosz, {\it Perturbations of Banach Algebras},
Lecture Notes in Mathematics {\bf 1120}, Springer-Verlag, 1985.

\bibitem {ja5} K.~Jarosz, {\it $H^\infty(D)$ is stable}, Journal of
the London Mathematical Society (2) {\bf 37} (1988), 490--498.

\bibitem {kt} Y.~Katznelson, {\it An Introduction to Harmonic
Analysis}, 3rd edition, Cambridge University Press, 2004.

\bibitem {km} J.~Kigami, {\it Analysis on Fractals}, Cambridge
University Press, 2001.

\bibitem {k-r-1} A.~Koranyi and H.~Reimann, {\it Quasiconformal
mappings on the Heisenberg group}, Inventiones Mathematicae {\bf 80}
(1985), 309--338.

\bibitem {k-r-2} A.~Koranyi and H.~Reimann, {\it Quasiconformal
mappings on CR manifolds}, in {\it Complex Geometry and Analysis},
59--75, Lecture Notes in Mathematics {\bf 1422}, 1990.

\bibitem {k-r-3} A.~Koranyi and H.~Reimann, {\it Foundations for the
theory of quasiconformal mappings on the Heisenberg group}, Advances
in Mathematics {\bf 111} (1995), 1--87.

\bibitem {kr} S.~Krantz, {\it Function Theory of Several Complex
Variables}, AMS Chelsea, 2001.

\bibitem {k-p-1} S.~Krantz and H.~Parks, {\it The Geometry of Domains
in Space}, Birkh\"auser, 1999.

\bibitem {k-p-2} S.~Krantz and H.~Parks, {\it A Primer of
Real-Analytic Functions}, 2nd edition, Birkh\"auser, 2002.

\bibitem {l} T.Laakso, {\it Ahlfors $Q$-regular spaces with arbitrary
$Q > 1$ admitting weak Poincar\'e inequality}, Geometric and
Functional Analysis {\bf 10} (2000), 111-123; erratum, {\bf 12}
(2002), 650.

\bibitem {l-v} O.~Lehto and K.~Virtanen, {\it Quasiconformal Mappings
in the Plane}, 2nd edition, Springer-Verlag, 1973.

\bibitem {m-s} R.~Mac\'{\i}as and C.~Segovia, {\it Lipschitz functions
on spaces of homogeneous type}, Advances in Mathematics {\bf 33}
(1979), 257--270.

\bibitem {mat} P.~Mattila, {\it Geometry of Sets and Measures in
Euclidean Spaces}, Cambridge University Press, 1995.

\bibitem {mo} R.~Montgomery, {\it A Tour of Subriemannian Geometries,
their Geodesics and Applications}, American Mathematical Society,
2002.

\bibitem {p-w} M.~Picardello and W.~Woess, editors, {\it Random Walks
and Discrete Potential Theory}, Cambridge University Press, 1999.

\bibitem {p} P.~Pansu, {\it M\'etriques de Carnot--Carath\'eodory et
quasiisom\'etries des espaces sym\'etriques de rang un}, Annals of
Mathematics (2) {\bf 129} (1989), 1--60.

\bibitem {r} C.~Rickart, {\it Natural Function Algebras},
Springer-Verlag, 1979.

\bibitem {ro1} R.~Rochberg, {\it Deformation Theory of Uniform
Algebras}, Proceedings of the London Mathematical Society (3) {\bf 39}
(1979), 93--118.

\bibitem {ro2} R.~Rochberg, {\it The disk algebra is rigid},
Proceedings of the London Mathematical Society (3) {\bf 39} (1979),
119--129.

\bibitem {ro3} R.~Rochberg, {\it Deformation theory for uniform
algebras: An introduction}, in {\it Proceedings of the Conference on
Banach Algebras and Several Complex Variables}, 209--216, American
Mathematical Society, 1984.

\bibitem {ro4} R.~Rochberg, {\it Deformation of uniform algebras on
Riemann surfaces}, Pacific Journal of Mathematics {\bf 121} (1986),
135--181.

\bibitem {ru1} W.~Rudin, {\it Fourier Analysis on Groups}, Wiley,
1990.

\bibitem {ru2} W.~Rudin, {\it Function Theory in the Unit Ball of
${\bf C}^n$}, Springer-Verlag, 1980.

\bibitem {se1} S.~Semmes, {\it Some Novel Types of Fractal Geometry},
Oxford University Press, 2001.

\bibitem {se2} S.~Semmes, {\it Mappings and spaces, 2},
math.CA/0612265.

\bibitem {sr1} J.-P.~Serre, {\it Cohomologie des groupes discrets}, in
{\it Prospects in Mathematics}, 77--169, Princeton University Press,
1971.

\bibitem {sr2} J.-P.~Serre, {\it Lie Algebras and Lie Groups}, 2nd
edition, Lecture Notes in Mathematics {\bf 1500}, Springer-Verlag,
2006.

\bibitem {si} G.~Simmons, {\it Introduction to Topology and Modern
Analysis}, Krieger, 1983.

\bibitem {st1} E.~Stein, {\it Singular Integrals and Differentiability
Properties of Functions}, Princeton University Press, 1970.

\bibitem {st2} E.~Stein, {\it Harmonic Analysis: Real-Variable
Methods, Orthogonality, and Oscillatory Integrals}, with the
assistance of T.~Murphy, Princeton University Press, 1993.

\bibitem {s-w} E.~Stein and G.~Weiss, {\it Introduction to Fourier
Analysis on Euclidean Spaces}, Princeton University Press, 1971.

\bibitem {stt} E.~Stout, {\it The Theory of Uniform Algebras}, Bogden
\& Quigley, 1971.

\bibitem {sz1} R.~Strichartz, {\it Analysis on fractals}, Notices of
the American Mathematical Society {\bf 46} (1999), 1199--1208.

\bibitem {sz2} R.~Strichartz, {\it Differential Equations on Fractals:
A Tutorial}, Princeton University Press, 2006.

\bibitem {t} M.~Taibleson, {\it Fourier Analysis on Local Fields},
Princeton University Press, 1975.

\bibitem {v-s-c} N.~Varopoulos, L.~Saloff-Coste, and T.~Coulhon, {\it
Analysis and Geometry on Groups}, Cambridge University Press, 1992.

\bibitem {w1} N.~Weaver, {\it Lipschitz algebras and derivations of
von Neumann algebras}, Journal of Functional Analysis {\bf 139}
(1996), 261--300.

\bibitem {w2} N.~Weaver, {\it Lipschitz Algebras}, World Scientific,
1999.

\bibitem {w3} N.~Weaver, {\it Lipschitz algebras and derivations, II:
Exterior Differentiation}, Journal of Functional Analysis {\bf 178}
(2000), 64--112; erratum, {\bf 186} (2001), 546.

\bibitem {w} W.~Woess, {\it Random Walks on Infinite Graphs and
Groups}, Cambridge University Press, 2000.





\end{thebibliography}
\end{document}